\documentclass[10pt,a4paper,reqno]{amsart}
\usepackage{amsfonts}
\usepackage{amsthm}
\usepackage{amsmath}
\usepackage{mathtools}
\usepackage{amscd}
\usepackage[utf8]{inputenc}
\usepackage{t1enc}
\usepackage[mathscr]{eucal}
\usepackage{indentfirst}
\usepackage{graphicx}
\usepackage{graphics}
\usepackage{esint}
\usepackage{pict2e}
\usepackage{epic}
\usepackage{float}
\usepackage{MnSymbol}
\usepackage{multirow}
\usepackage{shuffle} 
\usepackage[table,xcdraw]{xcolor} 
\usepackage{hyperref}
\usepackage{tcolorbox}
\numberwithin{equation}{section}
\usepackage[margin=2.9cm]{geometry}
\usepackage{enumitem}
\usepackage{epstopdf} 
\usepackage{bbm}
\usepackage{enumitem}

\usepackage[toc,page]{appendix}

\usepackage[noadjust]{cite}

\allowdisplaybreaks

\newmuskip\pFqmuskip
\newcommand*\pFq[6][8]{%
	\begingroup 
	\pFqmuskip=#1mu\relax
	\mathchardef\normalcomma=\mathcode`,
	\mathcode`\,=\string"8000
	\begingroup\lccode`\~=`\,
	\lowercase{\endgroup\let~}\pFqcomma
	{}_{#2}F_{#3}{\left[\genfrac..{0pt}{}{#4}{#5};#6\right]}%
	\endgroup
}
\newcommand{\pFqcomma}{{\normalcomma}\mskip\pFqmuskip}

\theoremstyle{plain}
\newtheorem{theorem}{Theorem}[section]
\newtheorem{lemma}[theorem]{Lemma}
\newtheorem{corollary}[theorem]{Corollary}
\newtheorem{proposition}[theorem]{Proposition}

\theoremstyle{definition}

\newtheorem{remark}[theorem]{Remark}

\newtheorem{example}[theorem]{Example}

\newcommand{\Snew}{\mathcal{S}^{\text{new}}}
\newcommand{\Smin}{\mathcal{S}^{\text{min}}}

\begin{document}
	
	\title{Linear relations between $L$-values of newforms and moments of elliptic $K$ integral}
	
	\author[Kam Cheong Au]{Kam Cheong Au}
	
	\address{University of Cologne \\ Department of Mathematics and Computer Science \\ Weyertal 86-90, 50931 Cologne, Germany} 
	
	\email{kau1@smail.uni-koeln.de}
	

	\begin{abstract} By systematically translating certain integrals involving moments of the elliptic integral into $L$-values of modular forms on $\Gamma_1(4), \Gamma(4)$ and $\Gamma_1(8)$, and then utilizing relations among the critical $L$-values of cuspidal newforms, we obtain a dimension formula for the space spanned by these integrals. 
	\end{abstract}
	
	\maketitle
	
	\markboth{Kam Cheong Au}{Linear relations between $L$-values of newforms and moments of elliptic $K$ integral}

	\section*{Introduction}

Critical $L$-values of modular forms on the upper half-plane encode deep arithmetic information~\cite{atkin1970hecke, damerell1970functions, shimura1976special, eichler1957verallgemeinerung, shimura1977periods, shimura1971introduction}. This naturally leads to the problem of studying the linear relations they satisfy. More concretely, for a congruence subgroup $\Gamma \subset \text{SL}_2(\mathbb{Z})$, an integer $k \geq 3$, and a number field $K$, consider the $K$-vector space
\begin{equation}\label{L_value_space}
	\Bigl\{ (2\pi)^{k-1-s} L(f,s) \,\Big|\, s \in \{1,\ldots,k-1\},\ f \ \text{a weight $k$ cusp form on }\Gamma \text{ with $q$-coefficients in }K \Bigr\}.
\end{equation}
Our aim is to investigate the dimension of this space (over $K$) and to understand the relations among its elements. \par

At present, mathematics offers little progress toward establishing meaningful \emph{lower} bounds for this space, as it reduces to difficult irrationality questions. In this work, we restrict ourselves to providing an \emph{upper} bound which, at least conjecturally, appears to be correct. We focus on three illustrative choices of congruence subgroups: $\Gamma(2) \cong \Gamma_1(4)$, $\Gamma(4)$, and $\Gamma_1(8)$, they serve to demonstrate that dimension counting in these simple examples is already not entirely trivial. \par

One motivation for studying the space~\eqref{L_value_space} is that its elements correspond to depth-one modular multiple zeta values (MZVs) on congruence subgroups~\cite{brown2014multiple, matthes2016elliptic, broedel2016relations}. In this sense, our problem generalizes the classical question of finding linear relations among MZVs~\cite{zhao2016multiple}. We hope our depth-one analysis may provide a foundation for investigating higher-depth modular MZVs. \par

A second motivation arises from the study of moments of elliptic integrals, one of the focus of this article. Let
$$K(m) := \int_0^1 \frac{1}{\sqrt{(1-t^2)(1-mt^2)}} \, \mathrm{d}t$$
denote the complete elliptic integral of the first kind, it can be either interpreted as a hypergeometric series or a weight $1$ modular function. Both perspectives endow integrals of the shape
\begin{equation}\label{intK_shape}
	\int_0^1 g(m)\, K(m)^{A_1} K(1-m)^{A_2} \,\mathrm{d}m, 
	\qquad A_1,A_2 \in \mathbb{N}, \quad g(m) \ \text{an algebraic function of $m$},
\end{equation}
with rich arithmetic structure. Such integrals appear in diverse contexts, including random walks~\cite{borwein2012densities,Wanthesis,borwein2013three}, scattering amplitudes~\cite{bailey2008elliptic}, and lattice sums~\cite{wan2016integrals, borwein2013lattice, rogers2015moments}. \par

The hypergeometric incarnation of $K(m)$ provides a powerful toolbox for analyzing integrals of the form~\eqref{intK_shape}; see, for instance,~\cite{wan2012moments, bailey2008elliptic, borwein2012densities, rogers2015moments, zhou2014legendre, zhou2018two}. On the other hand, the modular incarnation of $K(m)$ as a weight~1 form has received comparatively less attention. For certain algebraic functions $g(m)$, the corresponding integrals can be expressed as critical $L$-values of modular forms on congruence subgroups. While this connection was acknowledged in earlier works~\cite{rogers2015moments,wan2016integrals} and \cite[Chaps.~5--6]{borwein2013lattice}, there it remains a rather hollow terminology and its rich theory was not utilized. This article incorporates this theory into the picture, our main results might be difficult to obtain using classical hypergeometric techniques alone.

\begin{theorem}\label{theorem_dim}(Dimension formulas)
	Let $k\geq 3$ be an integer, \\
	(a) Let $\mathcal{I}_k(\Gamma_1(4))$ denote the $\mathbb{Q}$-vector space spanned by the real numbers $$\left\{\int_0^1 K(m)^{k-s-1} K(1-m)^{s-1} m^i  \:\mathrm{d}m \big| s\in \{1,\cdots,k-1\} , i\in \{0,\cdots,\lfloor \frac{k}{2}\rfloor - 2\}\right\},$$
	then $$\dim_\mathbb{Q} \mathcal{I}_k(\Gamma_1(4)) \leq \begin{cases} 2\lfloor \frac{k+2}{4} \rfloor,& k\equiv 0 \pmod{2} \\ \lfloor \frac{k}{2} \rfloor - 1, & k\equiv 1 \pmod{2}\end{cases}.$$
	(b) Abbreviate $X = m^{1/4}$, let $\mathcal{I}_k(\Gamma(4))$ denote the $\mathbb{Q}$-vector space spanned by the real numbers $$\left\{\int_0^1 K(m)^{k-s-1} K(1-m)^{s-1} \frac{X^i}{ X^3 (1 + X) (1 + X^2)} \:\mathrm{d}m \big| s\in \{1,\cdots,k-1\}, i\in \{0,\cdots,2k-2\}\right\},$$
	then $$\dim_\mathbb{Q} \mathcal{I}_k(\Gamma(4)) \leq k.$$
	(c) Abbreviate $X = \sqrt{1+\sqrt{1-m}}$, let $\mathcal{I}_k(\Gamma_1(8))$ denote the $\mathbb{Q}(\sqrt{2})$-vector space spanned by the real numbers $$\left\{\int_0^1 K(m)^{k-s-1} K(1-m)^{s-1} \frac{X^i}{(X-1) X^2 (X+1)^2 (X+\sqrt{2})} \:\mathrm{d}m \big| s\in \{1,\cdots,k-1\}, i\in \{0,\cdots,2k-2\}\right\},$$
	then $$\dim_{\mathbb{Q}(\sqrt{2})} \mathcal{I}_k(\Gamma_1(8)) \leq \begin{cases} 2k-1,& k\equiv 0 \pmod{2} \\ \lfloor \frac{3k}{2} \rfloor, & k\equiv 1 \pmod{2}\end{cases}.$$
\end{theorem}

We shall prove the three spaces above are actually the space (\ref{L_value_space}) adjoining certain Dirichlet $L$-values. We also expect these \textit{upper bounds in dimensions are tight}.\par
 
For $D$ a fundamental discriminant, let $\chi_D(n) := \left(\frac{D}{n}\right)$ to be Dirichlet character associated to field $\mathbb{Q}(\sqrt{D})$, $L_D(s)$ to be the Dirichlet $L$-function associated to the quadratic field $\mathbb{Q}(\sqrt{D})$, namely,
$$L_D(s) := L(\chi_D,s) = \sum_{n\geq 1} \left( \frac{D}{n}\right) n^{-s}.$$
The next theorem gives numbers in above three spaces that are "closed-form" constants, the list is complete in sense that all remaining basis element come from $L$-values of non-CM cusp forms.
\begin{theorem}\label{theorem_const}[Closed-form constants] Let $k\geq 3$ be an integer, denote
	\begin{equation}\label{Omega48}\Omega_{-4} := \frac{\Gamma \left(\frac{1}{4}\right)^2}{4 \pi ^{3/2}} ,\qquad \Omega_{-8} := \frac{\Gamma \left(\frac{1}{8}\right) \Gamma \left(\frac{3}{8}\right)}{2^{11/4} \pi ^{3/2}}.\end{equation}
	(a) The following numbers are in $\mathcal{I}_k(\Gamma_1(4))$:
	$$\begin{cases} \pi^{k-1}, \zeta(k-1), & \text{ when } k\equiv 0 \pmod{2} \\ \pi^{k-1} \Omega_{-4}^{k-1}, & \text{ when } k\equiv 1 \pmod{4}\end{cases}.$$
	(b) The following numbers are in $\mathcal{I}_k(\Gamma(4))$:
	$$\begin{cases} \pi^{k-1}, \zeta(k-1), & \text{ when } k\equiv 0 \pmod{2} \\ \pi^{k-1}, L_{-4}(k-1),\pi^{k-1} \Omega_{-4}^{k-1}, & \text{ when } k\equiv 1 \pmod{2}\end{cases}.$$
	(c) The following numbers are in $\mathcal{I}_k(\Gamma_1(8))$:
	$$\begin{cases} \pi^{k-1}, \zeta(k-1), L_8(k-1) & \text{ when } k\equiv 0 \pmod{2} \\ \pi^{k-1}, L_{-4}(k-1), L_{-8}(k-1), \pi^{k-1} \Omega_{-4}^{k-1}, \pi^{k-1} \Omega_{-8}^{k-1} & \text{ when } k\equiv 1 \pmod{4} \\
		\pi^{k-1}, L_{-4}(k-1), L_{-8}(k-1), \pi^{k-1} \Omega_{-8}^{k-1} & \text{ when } k\equiv 3\pmod{4}\end{cases}.$$
\end{theorem}
~\\[0.02in]

The organization of this article is as follows. Section~1 presents illustrative examples of the two main theorems to integrals (\ref{intK_shape}), recovering many isolated cases previously established in the literature. Section~2 begins the proof of our main theorems by first reducing them to statements about $L$-values. Section~3 addresses the contributions from Eisenstein series. Section~4 provides preliminaries on twisting cusp forms by Dirichlet characters. Finally, Section~5---the core of this work---analyzes the dimension of space (\ref{L_value_space}), our effort culminates in the proof of two theorems above.

\section*{Acknowledgment}
The author has received funding from the European Research Council (ERC) under the European Union’s Horizon 2020 research and innovation programme (grant agreement No. 101001179).

\section{Some numerical examples}
In this section we apply Theorems \ref{theorem_dim} and \ref{theorem_const} to small $k$, we will see many isolated results in literature on these integrals actually subsumes under them. \par
The proof of these two theorems can be made explicit: it gives explicit relations between elements in $\mathcal{I}_k(\Gamma)$, although not in any straightforward and easy-to-implement way. For illustration, we simply present such relations without proof.

We start with the simplest example.
\begin{example}
The case $k=4$ in part (a) of Theorem \ref{theorem_dim} gives $\dim_\mathbb{Q} \mathcal{I}_4(\Gamma_1(4)) \leq 2$, it contains the integrals
$$\int_0^1 K(m)^2 \:\mathrm{d}m,\quad \int_0^1 K(m) K(1-m) \:\mathrm{d}m.$$
Part (a) of Theorem \ref{theorem_const} says $\pi^3$ and $\zeta(3)$ lie in this space. Indeed, 
$$\int_0^1 K(m)^2 \:\mathrm{d}m = \frac{7\zeta(3)}{2},\quad \int_0^1 K(m) K(1-m) \:\mathrm{d}m= \frac{\pi^3}{8}.$$
These seem to be first mentioned in \cite{wan2012moments}, who obtained them using hypergeometric series. Since then many different proofs are known: see \cite{cantarini2019interplay, cantarini2022note}.  
\end{example}

\begin{example}
	The case $k=5$ in part (a) of Theorem \ref{theorem_dim} gives $\dim_\mathbb{Q} \mathcal{I}_5(\Gamma_1(4)) = 1$, combining with Theorem \ref{theorem_const}, we see each of
	$$\int_0^1 K(m)^i K(1-m)^{3-i} \:\mathrm{d}m ,\quad 0\leq i\leq 3$$
	is a rational multiple of $\pi^4 \Omega_{-4}^4= \dfrac{\Gamma(1/4)^8}{256 \pi ^2}$. Indeed, $$\int_0^1 K(m)^3 \:\mathrm{d}m = \frac{4}{5}\pi^4 \Omega_{-4}^4,\quad \int_0^1 K(m)^2 K(1-m) \:\mathrm{d}m = \frac{2}{3}\pi^4 \Omega_{-4}^4.$$
	These integrals have been proved in \cite{rogers2015moments} via manipulation of CM modular forms. It is also possible to prove them using purely hypergeometric techniques (\cite{zhou2014legendre}). 
\end{example}

\begin{example}\label{motiv_example}
	The case $k=6$ in part (a) of Theorem \ref{theorem_dim} gives $\dim_\mathbb{Q} \mathcal{I}_6(\Gamma_1(4)) \leq 4$, and Theorem \ref{theorem_const} says it contains $\zeta(5), \pi^5$. Explicitly, 
	$$\int_0^1 K(m)^4 (2m-1) \:\mathrm{d}m  = \frac{93}{8}\zeta(5),\quad \int_0^1 K(m)^3 K(1-m) (2m-1) \:\mathrm{d}m  = \frac{\pi^5}{64}.$$
	A hypergeometric proof of first one can be found in \cite{zhou2018two}. The other two constants comes from the $L$-values
	$$\int_0^1 K(m)^4 = 24L(f,5),\quad \int_0^1 K(m)^3 K(1-m) \:\mathrm{d}m  = 6\pi L(f,4),$$
	with $f(\tau) = \eta(2\tau)^{12} \in \Snew_6(\Gamma_0(4))$. Here we mention a curious representation of $L(f,5)$ in terms of hypergeometric series (see Appendix for a proof)
	$$L(f,5) =\frac{\pi^6}{1024} \: \pFq{9}{8}{\frac{1}{2},\frac{1}{2},\frac{1}{2},\frac{1}{2},\frac{1}{2},\frac{1}{2},\frac{1}{2},\frac{1}{2},\frac{5}{4}}{1,1,1,1,1,1,1,\frac{1}{4}}{1}.$$
\end{example}

\begin{example}
	The case $k=7$ in part (a) of Theorem \ref{theorem_dim} gives $\dim_\mathbb{Q} \mathcal{I}_7(\Gamma_1(4)) \leq 2$. So we only require two numbers to express each of $12$ integrals
	$$\int_0^1 x^i K(m)^j K(1-m)^{5-j} \:\mathrm{d}m,\quad 0\leq i\leq 1, 0\leq j\leq 5.$$
	Let $x_1 = \int_0^1 K(m)^5 \:\mathrm{d}m,\: x_2 = \int_0^1 K(m)^5 m \:\mathrm{d}m,$
	they are enough to express all cases, for example
	$$\int_0^1 K(m)^4 K(1-m) \:\mathrm{d}m = \frac{17}{30}x_1+\frac{11}{120}x_2, \quad \int_0^1 K(m)^3 K(1-m)^2 m \:\mathrm{d}m = -\frac{2}{15}x_1+\frac{31}{60}x_2.$$ 
Note that Theorem \ref{theorem_const} does not provide us with any "closed-form" constant in this case. The numbers $x_1, x_2$ are $L$-values of two cusp forms\footnote{Unlike example above, they are not Hecke eigenform, which is defined on $\mathbb{Q}(\sqrt{-15})$.}on $\Snew_7(\Gamma_1(4))$. 
\end{example}

\begin{example}
	The case $k=9$ in part (a) of Theorem \ref{theorem_dim} gives $\dim_\mathbb{Q} \mathcal{I}_9(\Gamma_1(4)) \leq 3$. So we only require three numbers, one being $\pi^8 \Omega_{-4}^8 =  \dfrac{\Gamma(1/4)^{16}}{2^{16} \pi ^4}$, to express each of $24$ integrals
	$$\int_0^1 x^i K(m)^j K(1-m)^{7-j} \:\mathrm{d}m ,\quad 0\leq i\leq 2, 0\leq j\leq 7.$$
	Let $x_1 = \int_0^1 K(m)^7 \:\mathrm{d}m, \: x_2 = \int_0^1 K(m)^7 m \:\mathrm{d}m,$
	they are enough to express all cases, for instances
	$$\int_0^1 K(m)^4 K(1-m)^3 m^2 \:\mathrm{d}m = \frac{4}{21}x_1+\frac{16}{35}x_2-\frac{432}{175} \pi^8 \Omega_{-4}^8.$$ 
\end{example}
~\\[0.02in]
Now we move on to examples for part (b) of Theorems \ref{theorem_dim} and \ref{theorem_const}.
\begin{example}
The case $k=3$ in part (b) of Theorem \ref{theorem_dim} gives $\dim_\mathbb{Q} \mathcal{I}_3(\Gamma(4)) \leq 3$, the three constants should be $\pi^2, L_{-4}(2), \pi^2 \Omega_{-4}^2$. For example, 
$$\int_0^1 K(m) \frac{m^{-3/4}}{(1+m^{1/2})(1+m^{1/4})} \:\mathrm{d}m = -\frac{\pi^2}{16}+L_{-4}(2)+\frac{\Gamma \left(\frac{1}{4}\right)^4}{16 \pi }.$$
\end{example}

\begin{example}
	The case $k=4$ gives $\dim_\mathbb{Q} \mathcal{I}_4(\Gamma(4)) \leq 4$, they are given by $\zeta(3),\pi^3$ and $x_1 = \int_0^1 K(m)^2 m^{-3/4}, \: x_2 = \int_0^1 K(m)^2 m^{-1/2}$. For example, 
	$$\int_0^1 K(1-m)^2 \frac{m^{-3/4}}{(1+m^{1/2})(1+m^{1/4})} \:\mathrm{d}m = \frac{7}{2}x_1+\frac{3}{4}x_2+\frac{\pi^3}{16}.$$
	$x_1, x_2$ are given by $L$-values of the unique newform $f(\tau) = \eta(2\tau)^4 \eta(4\tau)^4 \in \Snew_4(\Gamma_0(8))$, explicitly, $x_1 =  4L(f,3) + 4\pi L(f,2); x_2 = 8L(f,3)$.
\end{example}

\begin{example}
	The case $k=5$ gives $\dim_\mathbb{Q} \mathcal{I}_5(\Gamma(4)) \leq 5$, three of them are given by $\pi^4, L_{-4}(4)$ and $\pi^4 \Omega_{-4}^4$. Some special cases for which evaluation involves only these three numbers are 
	$$\begin{aligned}\int_0^1 K(m)^2 K(1-m) \frac{1}{(1+m^{1/4})(1+m^{1/2})} \:\mathrm{d}m &= \frac{\pi ^4}{320}+\frac{\Gamma \left(\frac{1}{4}\right)^8}{1280 \pi ^2}. \\
	\int_0^1 K(m)^3 m^{-1/2} \:\mathrm{d}m &= \frac{3 \Gamma \left(\frac{1}{4}\right)^8}{640 \pi ^2}. \\
	\int_0^1 K(m)^3 m^{1/2} \:\mathrm{d}m &= \frac{24 L_{-4}(4)}{5}+\frac{3 \Gamma \left(\frac{1}{4}\right)^8}{3200 \pi ^2}.\end{aligned}$$
	There are many variations of similar results in the literature (see \cite{zhou2014some, zhou2014legendre, rogers2015moments}), it seems they are all elements of space $\mathcal{I}_5(\Gamma(4))$ and so covered by our main theorems.
\end{example}

\begin{example}
	The case $k=7$ gives $\dim_\mathbb{Q} \mathcal{I}_7(\Gamma(4)) \leq 7$, three of them are given by $\pi^6, L_{-6}(6)$ and $\pi^6 \Omega_{-6}^6$, the remaining four constants comes from non-CM cusp forms. Among any five numbers in $\mathcal{I}_7(\Gamma(4))$, there exists a $\mathbb{Q}$-linear combination that eliminates these four constants. For example (with $X=m^{1/4}$):
	$$\begin{aligned}
	&\int_0^1 K(m)^5 \frac{201473 X^4+10177 X^3-14572 X^2-908 X+13664}{X^3(X+1)}  \:\mathrm{d}m = \frac{49 \pi ^6}{32}+\frac{1281 \Gamma \left(\frac{1}{4}\right)^{12}}{32 \pi ^3}. \\
	&\int_0^1 K(m)^2 K(1-m)^3 \frac{1037 - 4576 X + 1037 X^2 - 2197 X^3 + 14518 X^4}{X^3(X^2+1)}  \:\mathrm{d}m = \frac{1037 \Gamma \left(\frac{1}{4}\right)^{12}}{256 \pi ^3}-\frac{119 \pi ^6}{32}. \\
	&\int_0^1 K(m)^5\frac{15973982 X^4+6280350 X^3-718485 X^2-718485 X+745664 }{X^3 (X+1) \left(X^2+1\right)}  \:\mathrm{d}m \\ & \quad =  1283520 L_{-4}(6)-\frac{197629 \pi ^6}{64}+\frac{34953 \Gamma \left(\frac{1}{4}\right)^{12}}{16 \pi ^3}.
	\end{aligned}$$
	
\end{example}

Now we move on to examples for part (c) of Theorems \ref{theorem_dim} and \ref{theorem_const}.

\begin{example}
	The case $k=3$ gives $\dim_{\mathbb{Q}(\sqrt{2})} \mathcal{I}_3(\Gamma_1(8)) \leq 4$, they are given by $\pi^2, L_{-4}(2), L_{-8}(2)$ and $\pi^2 \Omega_{-8}^2$. Some instances (with $X=\sqrt{1+\sqrt{1-m}}$):
	$$\begin{aligned}\int_0^1 K(m) \frac{1}{(X-1) X^2 (X+1)^2 \left(X+\sqrt{2}\right)} \:\mathrm{d}m &= 4 (\sqrt{2}+1) L_{-4}(2)-\frac{8}{3} (\sqrt{2}+1) L_{-8}(2)-\frac{\sqrt{2}+1}{18} \pi ^2.\\
	\int_0^1 K(m) \frac{1}{(X - 1) X (1 + X)^2} \:\mathrm{d}m &= \frac{\pi ^2}{18}+\frac{\Gamma \left(\frac{1}{8}\right)^2 \Gamma \left(\frac{3}{8}\right)^2}{72 \pi }.\end{aligned}
	$$ 
\end{example}

\begin{example}
	The case $k=5$ gives $\dim_{\mathbb{Q}(\sqrt{2})} \mathcal{I}_5(\Gamma_1(8)) \leq 7$, five are them are given in Theorem \ref{theorem_const}. Among any three numbers in $\mathcal{I}_5(\Gamma_1(8)) $, there exists a linear combination that eliminates the remaining two constants. For example, there exists $c_i \in \mathbb{Q}(\sqrt{2})$ such that (with $X=\sqrt{1+\sqrt{1-m}}$)
	\begin{multline*}\int_0^1 \frac{c_1 K(m)^3 + c_2 K(m)^2 K(1-m) + c_3 K(m) K(1-m)^2}{X(X-1)(X+\sqrt{2})} \:\mathrm{d}m \\ = c_4 L_{-8}(4) + c_5 L_{-4}(4) + c_6\pi ^4 + c_7 \frac{\Gamma \left(\frac{1}{8}\right)^4 \Gamma \left(\frac{3}{8}\right)^4}{2048 \pi ^2} + c_8 \frac{\Gamma \left(\frac{1}{4}\right)^8}{256 \pi ^2}.\end{multline*}
	Explicitly, 
	{\small \allowdisplaybreaks \begin{alignat*}{2}c_1 &= 33174 \sqrt{2}+15922 \qquad &&c_2= -285 \sqrt{2}-33630, \\ c_3 &= -74469 \sqrt{2}-63927 \qquad &&c_4 = 128 (1327 \sqrt{2}-908), \\ c_5 &= \frac{912}{5} (454 \sqrt{2}-1327) \qquad &&c_6 = \frac{29}{4}-\frac{117}{8 \sqrt{2}}, \\ c_7 &= \frac{56}{9} (1477 \sqrt{2}+412) \qquad &&c_8 = \frac{2}{25} (1559027-1288929 \sqrt{2}).
	 \end{alignat*}}
\end{example}

	\section{Reduction to $L$-values}
	Main result of this section is Theorem \ref{int_to_L_values} that we prove at the end of this section, it connects $\mathcal{I}_k(\Gamma)$ to space spanned by $L$-values of modular forms. \par
	Throughout we write $q=e^{2\pi i \tau}$ for $\tau$ in upper half plane. Recall the three theta series $$\theta_2(\tau) = \sum_{n\in \mathbb{Z}} q^{(n+1/2)^2/2},\quad \theta_3(\tau) = \sum_{n\in \mathbb{Z}}q^{n^2/2},\quad \theta_4(\tau) = \sum_{n\in \mathbb{Z}} (-1)^n q^{n^2/2}.$$
	 
	 \begin{lemma}\label{theta234_in_m}
	 Let $\tau \in i\mathbb{R}^{>0}$, write $m = \dfrac{\theta_2^4(\tau)}{\theta_3^4(\tau)}$, so that $0<m<1$. Then \\ (a) $$K(m) = \frac{\pi}{2} \theta_3^2(\tau) ,\qquad iK(1-m) = \frac{\pi}{2} \tau \theta_3^2(\tau).$$ and $$\frac{d\tau}{dm} = \frac{\pi}{4i} \frac{1}{m(1-m)K(m)^2}.$$ 
	 (b) For $r\in  \{1,2,1/2\}$ and $i\in \{2,3,4\}$, we have the following expressions of $\frac{\pi}{2K(m)}\theta_i^2(r\tau)$ in terms of $m$. \begin{table}[H]\large
	 	\begin{tabular}{|c|c|c|c|}
	 		\hline
	 		$\frac{\pi}{2K(m)} \theta_i^2(r\tau)$ & $i=2$                   & $i=3$                    & $i=4$           \\ \hline
	 		$r=1$                                 & $\sqrt{m}$               & $1$                      & $\sqrt{1-m}$    \\ \hline
	 		$r=2$                                 & $\dfrac{1-\sqrt{1-m}}{2}$ & $\dfrac{1+\sqrt{1-m}}{2}$ & $\sqrt[4]{1-m}$ \\ \hline
	 		$r=\frac{1}{2}$                       & $2\sqrt[4]{m}$        & $1+\sqrt{m}$             & $1-\sqrt{m}$    \\ \hline
	 	\end{tabular}
	 	\label{tab:my-table}
	 \end{table}
	 \end{lemma}
	 \begin{proof}
	 Part (a) is a famous two-century-old result, see \cite{borwein1987pi,rogers2015moments,whittaker1920course}. For part (b), note that knowing two entries of a given row also gives the third entry because of the identity $\theta_2^4(\tau) + \theta_4^4(\tau) = \theta_3^4(\tau)$. The expression for $(r,i) = (1,2)$ and $(1,3)$ is contained in the statement of part (a). We can then obtain the entries for $(r,i) = (2,2)$ and $(2,3)$ by solving for $\theta_3^2(2\tau), \theta_2^2(2\tau)$ from identities
	 $$\theta_2^2(\tau) = 2\theta_3(2\tau) \theta_2(2\tau), \qquad \theta_2^4(2\tau) + \theta_3^4(2\tau) = \frac{1}{2}\theta_2^4(\tau) + \frac{1}{2}\theta_3^4(\tau),$$
	 and picking the solution which is positive on $0<m<1$. For the last row, replace $\tau$ by $\tau/2$ in above identities and then solve for $\theta_3^2(\tau/2), \theta_2^2(\tau/2)$. 
	 \end{proof}
	
	For a congruence subgroup $\Gamma$, we let $\mathcal{M}_k(\Gamma)$ denote the $\mathbb{C}$-vector space of weight $k$ modular form on $\Gamma$; for a field $K\subset \mathbb{C}$, denote $\mathcal{M}_k(\Gamma;K)$ to be forms whose $q$-coefficients are in $K$. For $\Gamma = \Gamma_1(N)$ or $\Gamma(N)$, we have
	 $$\mathcal{M}_k(\Gamma) = \mathcal{M}_k(\Gamma;\mathbb{Q}) \otimes_\mathbb{Q} \mathbb{C}$$
	 i.e. $\mathcal{M}_k(\Gamma)$ has a basis consists of rational Fourier expansions.  We quickly recall the structure of (holomorphic) modular forms on the congruence subgroups $\Gamma_1(4), \Gamma(4)$ and $\Gamma_1(8)$. 
	 \begin{lemma}
	 (a) $$\bigoplus_{k\geq 0} \mathcal{M}_k(\Gamma_1(4);\mathbb{Q})$$ is freely generated over $\mathbb{Q}$ by two generators $\theta_3^2(2\tau), \theta_4^4(2\tau).$ \\
	 (b) $$\bigoplus_{k\geq 0} \mathcal{M}_k(\Gamma(4);\mathbb{Q})$$ is generated over $\mathbb{Q}$ by three generators
	 	$\theta_3^2(2\tau), \theta_2^2(2\tau), \theta_2^2(\tau),$
	 	subject to a single weight $2$ relation.
	 (c) $$\bigoplus_{k\geq 0} \mathcal{M}_k(\Gamma_1(8);\mathbb{Q})$$ is generated over $\mathbb{Q}$ by three generators
	 $\theta_3^2(2\tau), \theta_4^2(2\tau), \theta_3(2\tau)\theta_3(4\tau),$ subject to a single weight $2$ relation.
	 \end{lemma}
	 \begin{proof}
	 It is straightforward to show these generators are indeed modular forms of corresponding subgroup. The dimension formulas
	 $$\dim(\mathcal{M}_k(\Gamma_1(4)) = \lfloor \frac{k}{2} \rfloor + 1,\qquad \dim(\mathcal{M}_k(\Gamma(4)) =  \dim(\mathcal{M}_k(\Gamma_1(8)) = 2k+1$$ then implies for $\Gamma_1(4)$, the two generators are free; for $\Gamma(4)$ and $\Gamma_1(8)$, they are related by a single relations at weight $2$.
	 \end{proof}

	 For a congruence subgroup $\Gamma$, let $\mathcal{M}^{(0,\infty)}_k(\Gamma, K)$ be the subspace of $\mathcal{M}_k(\Gamma, K)$ consisting of those that vanishes at cusps $0$ and $\infty$. For $f(\tau)$ defined on upper complex plane, we define $(B_2 f)(\tau)$ to be $f(\tau/2)$. 
	 \begin{lemma}Let $\tau \in i\mathbb{R}^{>0}$, write $m = \frac{\theta_2^4(\tau)}{\theta_3^4(\tau)}$. Also fix an integer $k\geq 3$, \\
	 (a) $B_2 \mathcal{M}_k(\Gamma_1(4);\mathbb{Q})$ has a $\mathbb{Q}$-basis given by
	 $$\left(\frac{2K(m)}{\pi}\right)^k m^i, \qquad i\in \{0,1,\cdots,\lfloor \frac{k}{2}\rfloor\}.$$ \\
	 (b) $B_2 \mathcal{M}_k(\Gamma(4);\mathbb{Q})$ has a $\mathbb{Q}$-basis given by
	 $$\left(\frac{2K(m)}{\pi}\right)^k m^{i/4}, \qquad i\in \{0,1,\cdots, 2k\}.$$ \\	 
	 (c) Let $K=\mathbb{Q}(\sqrt{2})$, $B_2 \mathcal{M}_k(\Gamma_1(8);K)$ has a $K$-basis given by
	 $$\left(\frac{2K(m)}{\pi}\right)^k (1+\sqrt{1-m})^{i/2}, \qquad i\in \{0,1,\cdots, 2k\}.$$
	 \end{lemma}
	 \begin{proof}
	 (a) From lemma above, $B_2\mathcal{M}_k(\Gamma_1(4);\mathbb{Q})$ is spanned by $$(\theta_3^2(\tau))^a (\theta_4^4(\tau))^b  = \left(\frac{2K(m)}{\pi}\right)^k (1-m)^b, \qquad a\geq 0, b\geq 0, a+2b=k.$$
	 where we used Lemma \ref{theta234_in_m}. This has the same span as $\left(\frac{2K(m)}{\pi}\right)^k m^i$. \\
	 (b) Abbreviate $X=m^{1/4}$, $B_2\mathcal{M}_k(\Gamma(4);\mathbb{Q})$ is spanned by $$(\theta_3^2(\tau))^a (\theta_2^2(\tau))^b  (\theta_2^2(\frac{\tau}{2}))^c = \left(\frac{2K(m)}{\pi}\right)^k (m^{1/2})^b (2m^{1/4})^c = \left(\frac{2K(m)}{\pi}\right)^k 2^c X^{c+2b},$$
	 for non-negative $a,b,c$ with $a+b+c=k$. Note that $c+2b \leq k+b\leq 2k$, so above are contained in the $(2k+1)$-dimensional space spanned by $\{(2K/\pi)^k X^i : 0\leq i \leq 2k\}$. The space $\mathcal{M}_k(\Gamma(4))$ has the same dimension, so this containment is an equality.\\
	 (c) Abbreviate $X=\sqrt{1+\sqrt{1+m}}$, $B_2\mathcal{M}_k(\Gamma_1(8));K)$ is spanned by $$(\theta_3^2(\tau))^a (\theta_4^2(\tau))^b  \left(\theta_3(2\tau)\theta_3(\tau)\right)^c = \left(\frac{2K(m)}{\pi}\right)^k (\sqrt{1-m})^b \left(\sqrt{\frac{1+\sqrt{1-m}}{2}}\right)^c = \left(\frac{2K(m)}{\pi}\right)^k 2^{-c/2} (X^2-1)^b X^c,$$
	 for non-negative $a,b,c$ with $a+b+c=k$. Again $c+2b \leq 2k$, so above are contained in the $(2k+1)$-dimensional space spanned by $\{(2K/\pi)^k X^i : 0\leq i \leq 2k\}$. The space $\mathcal{M}_k(\Gamma_1(8))$ has the same dimension, so this containment is an equality.
	 \end{proof}
	 
	 We introduce a notation that expedite our treatment for all three cases: let \begin{equation}\label{KandX}K= \begin{cases} \mathbb{Q} & \text{ for }\Gamma = \Gamma_1(4) \\  \mathbb{Q} & \text{ for }\Gamma = \Gamma(4) \\  \mathbb{Q}(\sqrt{2}) & \text{ for }\Gamma = \Gamma_1(8)    \end{cases}, \qquad X =  \begin{cases} m & \text{ for }\Gamma = \Gamma_1(4) \\  m^{1/4} & \text{ for }\Gamma = \Gamma(4) \\  \sqrt{1+\sqrt{1-m}} & \text{ for } \Gamma = \Gamma_1(8)    \end{cases}.\end{equation}
	 Also denote $V_n(K)$ be the $K$-vector space of univariate polynomials with degree $\leq n$, which has dimension $n+1$. 
	 \begin{lemma}\label{vec_iso}
	 For $\Gamma \in \{\Gamma_1(4),\Gamma(4),\Gamma_1(8)\}$, using the notations $K$ and $X$ as above. \\
	 (a) Let $n=\dim \mathcal{M}_k(\Gamma) - 1$, the map defined by $$\begin{aligned}P:\mathcal{M}_k(\Gamma;K) &\to V_n(K) \\ f(\tau) &\mapsto g, \qquad (B_2f)(\tau) := f(\tau/2) =  \left(\frac{2K(m)}{\pi}\right)^k g(X) \end{aligned}$$ is an isomorphism. \\
	 (b) Image of $\mathcal{M}_k^{(0,\infty)}(\Gamma;K)$ under $P$ consists of polynomial such that
	  $$\begin{cases} X(1-X) h(X) & \Gamma = \Gamma_1(4) \text{ or } \Gamma(4) \\  (X-1)(X-\sqrt{2}) h(X) & \Gamma = \Gamma_1(8)\end{cases},$$
	  here $h(X)$ is any polynomial in $V_{n-2}(K)$.
	 \end{lemma}
	 \begin{proof}
	 (a) is simply a reformulation of above lemma. The weight $1$ modular function $2K(m)/\pi = \theta_3(\tau)^2$ is non-vanishing at cusps $0$ and $\infty$, so $f$ vanishes at $0$ and $\infty$ iff the corresponding $g(X)$ vanishes at $m=0$ and $m=1$; for $X=m$ or $m^{1/4}$, this happens iff $g(X)$ is divisible by $X(1-X)$; for $X=\sqrt{1+\sqrt{1-m}}$, for which values of $m=0,1$ are $1$ and $\sqrt{2}$, this happens iff $g(X)$ is divisible by $(X-1)(X-\sqrt{2})$. 
	 \end{proof}

	 Now we claim the three spaces of integrals mentioned in Theorem \ref{theorem_dim} are equivalent to
	 \begin{equation}\label{aux1}\mathcal{I}_k(\Gamma) = \left\{ \int_0^1 K(m)^{k-s-1} K(1-m)^{s-1} \frac{g(X)}{m(1-m)} dm  \Big|  s\in \{1,\cdots,k-1\}, g\in P(\mathcal{M}_k^{(0,\infty)}(\Gamma;K)) \right\}.\end{equation}
	 Indeed, using notation of above lemma, $\Gamma = \Gamma_1(4)$, $\frac{g(X)}{m(1-m)} = \frac{X(1-X) h(X)}{X(1-X)} = h(X)$; for $\Gamma = \Gamma(4)$, $$\frac{g(X)}{m(1-m)} = \frac{X(1-X) h (X)}{X^4(1-X^4)} = \frac{h(X)}{X^3(1+X)(1+X^2)};$$
	 for $\Gamma = \Gamma_1(8)$, $$\frac{g(X)}{m(1-m)} = \frac{(X-1)(X-\sqrt{2}) h(X)}{(2 X^2 - X^4) (1 - 2 X^2 + X^4)} = \frac{h(X)}{(X-1) X^2 (X+1)^2 \left(X+\sqrt{2}\right)};$$
~\\[0.02in]
	 For $f = \sum_{n\geq 0} a_n q^n \in \mathcal{M}^{(0,\infty)}_k(\Gamma)$ with $\Gamma$ some congruence subgroup, we define its $L$-function to be
	 $$L(f,s) = \frac{(2\pi)^s}{\Gamma(s)} \int_0^\infty f(ix) x^s \frac{dx}{x}.$$
	 This is an entire function in $s$. When $f(\tau) = \sum_{n\geq 0} a_n(f) q^n$ has $q$-expansion consists of integral powers of $q$, this is the same as the familiar $L(f,s) = \sum_{n\geq 1} a_n(f) n^{-s}$. \par
	 
	 From $$L(f,s) = \frac{\pi^s}{\Gamma(s)} \int_0^\infty (B_2 f)(ix) x^s \frac{dx}{x},$$
	 make the substitution $x = K(1-m)/K(m)$, then from part (a) of Lemma \ref{theta234_in_m}, we have
	 $$\begin{aligned}
	L(f,s) &= \frac{\pi^s}{\Gamma(s)} \frac{\pi}{4} \int_0^1 \frac{K(1-m)^{s-1}}{K(m)^{s+1}} \frac{(B_2f)(ix)}{m(1-m)} dm \\
	&= \frac{2^{k-2}\pi^{s+1-k}}{\Gamma(s)}\int_0^1 K(m)^{k-s-1} K(1-m)^{s-1} \frac{(P\circ f)(X)}{m(1-m)} dm
\end{aligned}$$
with $P$ the map in Lemma \ref{vec_iso}. Comparing this with equation (\ref{aux1}), we established
\begin{theorem}\label{int_to_L_values}
	 $$\mathcal{I}_k(\Gamma) = \text{Span}_K \left\{(2\pi)^{k-1-s} L(f,s) \Big| s\in\{1,\cdots,k-1\}, f\in \mathcal{M}_k^{(0,\infty)}(\Gamma;K)\right\},$$
	where the field $K$ is given by equation (\ref{KandX}). \end{theorem}

	 We continue the investigation by splitting space of modular forms into Eisenstein and cuspidal part, 
	 $$\mathcal{M}_k^{(0,\infty)}(\Gamma;K) = \mathcal{E}_k^{(0,\infty)}(\Gamma;K) \oplus  \mathcal{S}_k(\Gamma;K),$$
	 where $\mathcal{E}_k^{(0,\infty)}(\Gamma;K) = \mathcal{E}_k(\Gamma;K) \cap \mathcal{M}_k^{(0,\infty)}(\Gamma;K)$.
	 In next two sections shall separately analyze the contribution of $L$-values when $f$ lies in each part.
	 \\[0.02in] 
	 In remainder of this section, we assemble some well-known notations and facts. For a function defined on upper-half plane, recall its weight $k$ $\text{GL}_2(\mathbb{Q})$-action:
	$$f|_k \gamma = (\det \gamma)^{k/2} (c\tau +d)^{-k} f(\frac{a\tau+b}{c\tau+d}).$$
	 For a positive integer $N$ and $\chi$ a Dirichlet character modulo $N$, recall the space $$\mathcal{M}_k(N,\chi) := \left\{f: f|_k \gamma  = \chi(d) f, \gamma = \begin{psmallmatrix}a&b\\ c&d \end{psmallmatrix} \in \Gamma_0(N)\right\}.$$
	 If $\chi$ is a character modulus some divisor of $N$, then we still denote $\mathcal{M}_k(N,\chi)$ by extending $\chi$ modulo $N$. We denote $\mathbbm{1}$ as the trivial character (of some modulus clear in context). Also recall the following elementary fact (\cite[p.~260]{cohen2017modular}), for any positive integer $N$, the map \begin{equation}\label{GammaNinflation}\mathcal{M}_k(\Gamma(N)) \to \bigoplus_{\chi \mod N} \mathcal{M}_k(N^2,\chi), \qquad f(\tau) \mapsto f(N\tau)\end{equation}
	 is an isomorphism, where $\chi$ ranges over all characters modulo $N$. This implies $\mathcal{M}_k(\Gamma(4))$ is isomorphic to $\mathcal{M}_k(16,\mathbbm{1})$ for $k$ even and $\mathcal{M}_k(16,\chi_{-4})$ for $k$ odd.
	 
	 \section{$L$-values from Eisenstein subspace}
	 Recall our notation of the field $K$ as in equation (\ref{KandX}). Goal of this section is to prove the following theorem.
	 	 \begin{theorem}\label{eis_const_theorem}For $\Gamma = \Gamma_1(4), \Gamma(4)$ or $\Gamma_1(8)$, 
	 	\begin{equation}\label{eis_span_eqn}\text{Span}_K \left\{(2\pi)^{k-1-s} L(f,s) \Big| s\in\{1,\cdots,k-1\}, f\in \mathcal{E}_k^{(0,\infty)}(\Gamma;K) \right\}\end{equation} is contained in $K$-span of following numbers
	 	\begin{equation}\label{eis_consts}\begin{cases}\{\pi^{k-1},\zeta(k-1)\}& k\text{ even and }\Gamma = \Gamma_1(4) \\  \{ \}& k\text{ odd and }\Gamma = \Gamma_1(4)\\ \{\pi^{k-1},\zeta(k-1)\} & k\text{ even and }\Gamma = \Gamma(4) \\  \{\pi^{k-1},L_{-4}(k-1)\} & k\text{ odd and }\Gamma = \Gamma(4) \\ \{\pi^{k-1},\zeta(k-1),L_8(k-1)\} & k\text{ even and }\Gamma = \Gamma_1(8) \\ \{\pi^{k-1},L_{-4}(k-1),L_{-8}(k-1)\} & k\text{ odd and }\Gamma = \Gamma_1(8)  \end{cases}.\end{equation}
	 \end{theorem}
	 
	 We can quickly dispense the case $k$ odd and $\Gamma = \Gamma_1(4)$: $\Gamma_1(4)$ has three $\text{SL}_2(\mathbb{Z})$-equivalent cusps $\{0,1,1/2\}$, the cusp $1/2$ being irregular. Modular forms of odd weight vanish at irregular cusps, thus $\mathcal{E}_k^{(0,\infty)}(\Gamma,K)$ is trivial, since any non-trivial element would be cuspidal. \par
	 
	 For remaining five cases, the dimension of $\mathcal{E}_k(\Gamma)$ equals the number of $\text{SL}_2(\mathbb{Z})$-inequivalent cusps, so $\dim \mathcal{E}_k^{(0,\infty)}(\Gamma)$ equals this number minus $2$. Let $\psi$ be a primitive character modulo $u$, $\varphi$ a primitive character modulo $v$ such that $(\psi \varphi)(-1) = (-1)^k$, there exists (\cite[chap.~7]{miyake2006modular},\cite[chap.~4]{diamond2005first}) an Eisenstein series $E_k^{\psi,\varphi}(\tau) \in \mathcal{E}_k(uv, \psi \varphi)$ such that its $L$-function is $L(\psi,s) L(\varphi,s-k+1)$. For a positive integer $t$, denote $E_k^{\psi,\varphi,t} (\tau) := E_k^{\psi,\varphi}(t\tau)$. Then a basis for our Eisenstein spaces of concern is given by
	 \begin{equation}\label{eis_basis}\begin{cases}E_k^{\mathbbm{1},\mathbbm{1},1}, E_k^{\mathbbm{1},\mathbbm{1},2}, E_k^{\mathbbm{1},\mathbbm{1},4}& \mathcal{E}_k(\Gamma_1(4)), k\text{ even } \\  E_k^{\mathbbm{1},\mathbbm{1},1}, E_k^{\mathbbm{1},\mathbbm{1},2} , E_k^{\mathbbm{1},\mathbbm{1},4} ,E_k^{\mathbbm{1},\mathbbm{1},8} ,E_k^{\mathbbm{1},\mathbbm{1},16}, E_k^{\chi_{-4},\chi_{-4},1}   & \mathcal{E}_k(16,\mathbbm{1}), k\text{ even } \\E_k^{\chi_{-4},\mathbbm{1},1}, E_k^{\chi_{-4},\mathbbm{1},2} , E_k^{\chi_{-4},\mathbbm{1},4}, E_k^{\mathbbm{1},\chi_{-4},1}, E_k^{\mathbbm{1},\chi_{-4},2}, E_k^{\mathbbm{1},\chi_{-4},4}& \mathcal{E}_k(16,\chi_{-4}), k\text{ odd }  \\
	 E_k^{\mathbbm{1},\mathbbm{1},1},E_k^{\mathbbm{1},\mathbbm{1},2},E_k^{\mathbbm{1},\mathbbm{1},4},E_k^{\mathbbm{1},\mathbbm{1},8}, E_k^{\mathbbm{1},\chi_8,1},E_k^{\chi_8,\mathbbm{1},1}& \mathcal{E}_k(\Gamma_1(8)), k\text{ even } \\  E_k^{\chi_{-4},\mathbbm{1},1}, E_k^{\chi_{-4},\mathbbm{1},2}, E_k^{\mathbbm{1},\chi_{-4},1}, E_k^{\mathbbm{1},\chi_{-4},2}, E_k^{\mathbbm{1},\chi_{-8},1}, E_k^{\chi_{-8}, \mathbbm{1},1}   &\mathcal{E}_k(\Gamma_1(8)), k\text{ odd } \end{cases}.\end{equation}
	 
	 Recall the definition of Fricke involution on $\mathcal{M}_k(\Gamma_1(N))$: $$(W_N f)(\tau) := i^k N^{-k/2} \tau^{-k} f(-\frac{1}{N\tau}) = i^k f|_k \begin{psmallmatrix}&-1\\ N& \end{psmallmatrix},$$ and the completed $L$-function $$\Lambda_N(f,s) := N^{s/2} (2\pi)^{-s} \Gamma(s) L(f,s),$$
	 then $\Lambda_N(\omega_N f, s) = \Lambda_N(f,k-s)$ for any $f\in \mathcal{M}_k(\Gamma_1(N))$ (\cite[chap.~5]{diamond2005first}).
	 
	 We illustrate the proof of the theorem in the first case. \begin{proof}[Proof of Theorem \ref{eis_const_theorem} for $\mathcal{E}_k(\Gamma_1(4))$ and $k$ even]
	  In this case, $E_k^{\mathbbm{1},\mathbbm{1}}$ is simply the familiar $\text{SL}_2(\mathbb{Z})$ Eisenstein series, by consider the behaviour under transformation $\tau\mapsto -1/\tau$, one readily sees $\mathcal{E}_k^{(0,\infty)}(\Gamma)$ is spanned by single element $$f = E_k^{\mathbbm{1},\mathbbm{1},1} - (1+2^k) E_k^{\mathbbm{1},\mathbbm{1},2} + 2^k E_k^{\mathbbm{1},\mathbbm{1},4},$$
	 thus the equation (\ref{eis_span_eqn}) becomes
	 $$\text{Span}_K \left\{(2\pi)^{k-1-s} \left(1-(1+2^k)2^{-s} + 2^{k} 4^{-s}\right)\zeta(s)\zeta(s-k+1) \Big| s\in\{1,\cdots,k-1\} \right\}.$$
	 One checks for $s=k-1$, above is $\mathbb{Q}\zeta(k-1)$; for $s=k-2,k-4,\cdots,2$, it gives $\mathbb{Q}\pi^{k-1}$; for $s=k-3,k-5,\cdots,3$, above is $0$ because $\zeta(-2\mathbb{N}) = 0$; the final case $s=1$ is slightly non-obvious, due to pole of $\zeta(s)$ at $s=1$. First method is simply taking limit $s\to 1$ and gives its value $\mathbb{Q} \pi^{k-2} \zeta'(2-k) = \mathbb{Q} \zeta(k-1)$. Second method, which actually generalize to all our five cases, is to rewrite equation (\ref{eis_span_eqn}) as
	 $$\text{Span}_K \left\{ (2\pi)^{k-1} \Lambda_4(f,s) \Big| s\in\{1,\cdots,k-1\} \right\}.$$
	 Since $W_N$ maps $\mathcal{E}_k^{(0,\infty)}(\Gamma)$ to itself, $W_N f = \pm f$, hence the functional equation $\Lambda_4(f,s) = \Lambda_4(f,k-s)$ says the case $s=1$ is equivalent to $s=k-1$. Hence its value at $s=1$ is in $\mathbb{Q}\zeta(k-1)$.
	 \end{proof}
	 
	 In order to handle all cases, we require the following lemma
	 \begin{lemma}
	 Let $k\geq 3$, $\psi$ be a primitive character modulo $u$, $\varphi$ a primitive character modulo $v$ such that $(\psi \varphi)(-1) = (-1)^k$. Denote $g(\psi) = \sum_{l=0}^{u-1} \psi(l) e^{2\pi i l/u}$ be the Gauss sum. (a) We have
	 $$W_{uv} (E_{k}^{\psi,\varphi}) =(\frac{u}{v})^{(1-k)/2}  \left( i^{-k} \frac{g(\psi)g(\varphi)}{\sqrt{uv}} \right)E_k^{\overline{\varphi},\overline{\psi}},$$
	 	(b) If $t,N\in \mathbb{N}$ such that $tuv \mid N$, then $$W_N E_k^{\psi,\varphi,t} =(-1)^k \sqrt{\frac{u}{v}} \frac{N^{k/2}}{(ut)^k} \left( i^{-k} \frac{g(\psi)g(\varphi)}{\sqrt{uv}} \right) E_k^{\overline{\varphi},\overline{\psi},N/(uvt)}.$$
	 \end{lemma}
	 \begin{proof}
	 (b) easily follows from (a), whose proof can be found in \cite[p.~3]{weisinger1977some}.
	 \end{proof}
	 
	 \begin{proof}[Proof of Theorem \ref{eis_const_theorem}]
	 Let $K$ be the field given in equation \ref{KandX} which depends on $\Gamma$. For real characters $\psi$ and $\varphi$, the number $ i^{-k} \frac{g(\psi)g(\varphi)}{\sqrt{uv}}$ in above lemma is either $1$ or $-1$, because $W_{uv}$ is an involution. Let $\mathcal{E}_k$ denote $K$-span of one of Eisenstein spaces given in equation (\ref{eis_basis}). \par
	 \underline{Claim 1: $W_N$ acts $\mathcal{E}_k$.}
	 We have to show the number $\sqrt{\frac{u}{v}} \frac{N^{k/2}}{(ut)^k}$ in part (b) of lemma above is always in $K$. As examples, for $\mathcal{E}_k(16,\mathbbm{1})$, either $(u,v,N) = (1,1,16)$ or $(u,v,N) = (4,4,16)$, so this always in $\mathbb{Q} = K$; for $\mathcal{E}_k(\Gamma_1(8))$, $(u,v,N) = (1,1,8), (1,8,8)$ or $(8,1,8)$, the number is always in $\mathbb{Q}(\sqrt{2}) = K$. The other cases are done similarly. \par
	 \underline{Claim 2: $W_N$ acts on $\mathcal{E}_k^{(0,\infty)}$.}
	 An $f\in \mathcal{E}_k$ is in this subspace iff $\lim_{\Im(\tau)\to \infty} f(\tau) = \lim_{\Im(\tau)\to \infty} W_Nf(\tau) = 0.$ Recall that (\cite[chap.~4]{diamond2005first}), $$\lim_{\Im(\tau)\to \infty} E_k^{\psi,\varphi}(\tau) = \begin{cases}L(1-k,\varphi)/2\qquad &\psi = \mathbbm{1} \\  0\qquad &\psi \neq \mathbbm{1}\end{cases}.$$
	 So $\mathcal{E}_k^{(0,\infty)}$ (\textit{a priori} only a complex space) is defined as kernel of some linear equation with coefficient in $K$.\par
	 \underline{Claim 3: $\mathcal{E}_k^{(0,\infty)}$ has a basis that are eigenvector of $W_N$.} This is because the eigenspaces are given by $\ker(W_N \pm 1)$, which are operators defined over $K$.\par
Hence the space (\ref{eis_span_eqn}) coincides with
$$\text{Span}_K \left\{ (2\pi)^{k-1} \Lambda_N(f,s) \Big| s\in\{1,\cdots,k-1\}, W_N f = \pm f \right\}.$$
Now the problem becomes easy, we do as an example the case $\mathcal{E}_k(\Gamma_1(8))$ with $k$ even. In this case, $(2\pi)^{k-1}\Lambda_8(f,s) $ equals $$ (2\pi)^{k-1-s}8^{s/2}\Gamma(s) \left[ A_1(s) \zeta(s)\zeta(s-k+1) + A_2(s) \zeta(s)L_8(s-k+1) + A_3(s) L_8(s)\zeta(s-k+1) \right]$$
where $A_i(s)$ are some Euler factors that is in $K$ when $s\in \mathbb{Z}$ (by Claim 3 above). Plugging in $s=k-1,k-2,\cdots,2$, one checks their values are in $K$-span of $\{\pi^{k-1},\zeta(k-1),L_8(k-1)\}$. For $s=1$, the pole of $\zeta(s)$ gives problem, but $W_N f =\pm f$ implies the functional equation $\Lambda_8(f,s) = \pm \Lambda_8(f,k-s)$, so it gets reduced to the case $s=k-1$. The other cases follow by similar examinations.
	 \end{proof}

	 \section{Preliminary on twist of newforms}
	 Here we recall some facts about twisting newforms by Dirichlet character. Let $N$ be a positive integer. For a newform (i.e. normalized Hecke eigenform) $f\in \Snew_k(N,\chi)$, a primitive Dirichlet character $\psi$ modulo $M$ (might or might not be the same as $N$), there exists a newform $f\otimes \psi \in \Snew_k(N', \chi \psi^2)$ on some level $N' \mid NM^2$ such that $$a_p(f\otimes \psi) = \psi(p) a_p(f)$$ for all primes $p$ not divisible by $NM^2$. By strong multiplicity one (\cite[306]{bump1998automorphic}), this $f\otimes \psi$ is uniquely determined. We say a newform $f\in \Snew_k(N,\psi)$ is twist-minimal if it is not a twist of any newform of lower level, the subspace spanned by twist-minimal newforms is denoted by $\Smin_k(N,\psi)$. 
	 
	 \begin{proposition}\label{twist_min_level}\cite[Theorem~3]{child2022twist}
	 Let $f\in \Smin_k(N,\chi)$ be a twist-minimal form, $\psi$ be a primitive Dirichlet character, and let $\psi'$ be the primitive character that induces $\chi \psi$. Then $f\otimes \psi \in \Snew_k(M,\chi\psi^2)$ with $$M = \text{LCM}(N,\text{cond}(\psi) \text{cond}(\psi')).$$
	 \end{proposition}

	 For our application, we need the following cases when $N$ divides $16$. In all what follows, we shall freely use dimension formulas for cusp space and newspace as found in \cite[chap.~6]{stein2007modular}. 
	 \begin{lemma}\label{prop_level16_twists}
	 Let $k\geq 3$ be an integer. \\
	 (a) $\Snew_k(N,\mathbbm{1}) = \Smin_k(N,\mathbbm{1})$ for $N=1,2,4,8$; $\Smin_k(16,\mathbbm{1}) = 0.$ \\
	 (b) $\Snew_k(4,\chi_{-4}) = \Smin_k(4,\chi_{-4})$, this space is closed under twisting by $\chi_{-4}$.\\
	 (c) $\Snew_k(8,\chi_{\pm 8}) = \Smin_k(8,\chi_{\pm 8})$, this space is closed under twisting by $\chi_{\pm 8}$.\\
	 (d) $\Snew_k(16,\chi_{-4}) = \Smin_k(16,\chi_{-4})$, this space is closed under twisting by $\chi_{-4}$.
	 \end{lemma}
	 \begin{proof}
	 (a) For $N=1,2,4,8$, this follows from from Proposition \ref{twist_min_level} by checking that twisting a form of lower level will never yield newform of that level with trivial character. For $N=16$, Proposition \ref{twist_min_level} again implies, for $d=1,2,4$ or $8$, twisting by $\chi_{-4}$ on $\Snew_k(d,\mathbbm{1})$ lands in $\Snew_k(16,\mathbbm{1})$. The equality of dimension $$\dim \Snew_k(16,\mathbbm{1}) = \sum_{d\mid 8}  \dim \Snew_k(d,\mathbbm{1})$$ says every newform on $\Snew_k(16,\mathbbm{1}) $ arises in this way, so there is no twist-minimal form. \\ 
	 (b) and (c) follows directly from Proposition \ref{twist_min_level}. \\
	 (d) It suffices to establish $\Snew_k(16,\chi_{-4}) = \Smin_k(16,\chi_{-4})$, the second assertion then follows from Proposition \ref{twist_min_level}. If this is not the case, say $f$ is not twist-minimal in $\Snew_k(16,\chi_{-4})$, then $f\otimes \chi_{-4} \in \Snew_k(16,\chi_{-4})$ for a newform $f\in \Snew_k(4,\chi_{-4})$ or $\Snew_k(8,\chi_{-4})$. The former case is excluded by (b) and the latter case is also not possible since $\Snew_k(8,\chi_{-4}) = \{0\}$ by dimension formula. 
	 \end{proof}
	 
	 A newform $f$ for which there exists a non-trivial character $\chi$ such that $f\otimes \chi = f$ occupies a special place in our investigation. When weight $k\geq 2$, this $\chi$ is then character of an imaginary quadratic field, $f$ is said to have complex multiplication (CM), the CM discriminant of $f$ is the discriminant of this field. It is easy to characterize CM newforms of level divisible by $16$.
	 \begin{proposition}\label{CM_forms}Let $k\geq 3$ be a positive integer. \\
	 	(a) $\Snew_k(\Gamma_1(4))$ contains a CM newform only when $k\equiv 1\pmod{4}$, this form is unique, has CM discriminant $-4$, and equals
	 	\begin{equation}\label{CM_level4}\frac{1}{4}\sum_{x,y\in \mathbb{Z}^2} (x+\sqrt{-1}y)^{k-1} q^{x^2+y^2} \in \Snew_k(4,\chi_{-4}).\end{equation}
	 	(b) $\Snew_k(\Gamma_1(8))$ contains a CM newform only when $k\equiv 1\pmod{2}$, this form is unique, has CM discriminant $-8$, and equals
	 	\begin{equation}\label{CM_level8}\frac{1}{2}\sum_{x,y\in \mathbb{Z}^2} (x+\sqrt{-2}y)^{k-1} q^{x^2+2y^2} \in \Snew_k(8,\chi_{-8}).\end{equation}
	 	(c) $\Snew_k(\Gamma_1(16))$ contains a CM newform only when $k\equiv 3\pmod{4}$, this form is unique, has CM discriminant $-4$, and equals
	 	\begin{equation}\label{CM_level16}\frac{1}{2}\sum_{\substack{x,y\in \mathbb{Z}^2 \\ x\text{ odd},\, y\text{ even }}} (x+\sqrt{-1}y)^{k-1} q^{x^2+y^2} \in \Snew_k(16,\chi_{-4}).\end{equation}
	 \end{proposition}
	 \begin{proof}
	 This follows from characterization of CM form in terms of certain Hecke character of imaginary quadratic field. See 
	 \cite[chap.~4.8]{miyake2006modular}, \cite[chap.~2.3]{cohen2017modular}, \cite{li2018computing}.
	 \end{proof}
	 
	 Recall our definition of numbers $\Omega_{-4}$ and $\Omega_{-8}$ as in equation (\ref{Omega48}).
	 \begin{proposition}\label{CM_L_value}
	 Let $f$ be one of CM newform listed in above Proposition, then
	 $$L(f,s) \in \mathbb{Q} \pi^s |D|^{(s-1)/2} \Omega_D^{k-1}, \qquad s\in \{1,\cdots,k-1\},$$
	 where $D$ is the corresponding CM discriminant.
	 \end{proposition}
	 \begin{proof}
	 This is not difficult from above explicit expressions (\ref{CM_level4}), (\ref{CM_level8}), (\ref{CM_level16}), then $L(f,k-1)$ becomes simply value of some Eisenstein series at a CM point (see \cite[p.~89]{bruinier20081} for details). Values at other $s$ then follows from Proposition \ref{Shimura_L_ratio} below.\end{proof}
	 
	 \begin{remark}
	 More generally, $L$-function of CM form is equivalent to that of a certain Hecke character, and its critical $L$-value are expressible in terms of a (proven) transcendental number known as Chowla-Selberg period. For a negative fundamental discriminant, following \cite[p.~84]{bruinier20081}, denote $$\Omega_D := \frac{1}{\sqrt{2\pi |D|}} \left(\prod_{j=1}^{|D|-1} \Gamma(\frac{j}{|D|})^{\chi_D(j)} \right)^{w_D/4h_D}$$
	 with $h_D$ the class number and $w_D$ the number of units in $\mathbb{Q}(\sqrt{D})$. The numbers $\Omega_{-4}$ and $\Omega_{-8}$ given above are special cases of above formula. For any CM newform of weight $k$ with CM discriminant $D$, one has (\cite{damerell1970functions})  $$L(f,s) \in \overline{\mathbb{Q}} \pi^s \Omega_D^{k-1}.$$\end{remark}
	 ~\\[0.02in]
	 We are interested in twist of newform because it gives relations between $L$-values. For a primitive character $\chi$ modulo $N$, denote its Gauss sum
	 $$g(\chi) := \sum_{r=0}^{N-1} \chi(r) e^{2\pi i r/N}.$$
	 For a non-primitive character $\chi$, $g(\chi)$ is defined as Gauss sum of its primitive part.
	 For an automorphism $\sigma \in  \text{Aut}(\mathbb{C}/\mathbb{Q})$ and $f\in \mathcal{M}_k(\Gamma)$, we write $f^\sigma \in \mathcal{M}_k(\Gamma)$ to be the form with $q$-coefficient those of $f$ acted by $\sigma$. Paramount in our investigation if the following due to Shimura. 
	 \begin{proposition}\cite[Theorem~1]{shimura1976special}\label{Shimura_L_ratio}
	 Let $k\geq 3$ be an integer, $f$ be a newform of some level, $\chi$ a Dirichlet character of some modulus. Define $$A(f,\chi,m) := \frac{(2\pi i)^{k-1-m}}{g(\chi)} L(f\otimes \chi, m).$$
	 Then for any $\sigma \in \text{Aut}(\mathbb{C}/\mathbb{Q})$; integers $m_1,m_2$ with $1\leq m_1,m_2\leq k-1$; $(\chi_1 \chi_2)(-1) = (-1)^{m_1 - m_2}$, we have, whenever denominator does not vanish,
	 $$\left(\frac{A(f,\chi_1,m_1)}{A(f,\chi_2,m_2)}\right)^\sigma =\frac{A(f^\sigma,\chi_1^\sigma,m_1)}{A(f^\sigma,\chi_2^\sigma,m_2)}.$$
	 In particular, the ratio $\frac{A(f,\chi_1,m_1)}{A(f,\chi_2,m_2)}$ is algebraic and lies in the field $\mathbb{Q}(f,\chi_1,\chi_2)$. 
	 \end{proposition} 
	 
	 \begin{example}\label{k-1andk-2reduction}
	 Take $\chi_1=\chi_2 = \mathbbm{1}$, then above Proposition implies the well-known fact
	 $$\frac{(2\pi)^{m_2-m_1} L(f,m_1)}{L(f,m_2)} \in \mathbb{Q}(f)$$
	 whenever $1\leq m_1,m_2\leq k-1$ (i.e. they are critical $L$-values) and $m_1,m_2$ have same parity. Thus up to $\overline{\mathbb{Q}}$-multiple, only two numbers $L(f,k-1)$ and $(2\pi) L(f,k-2)$ appear among the $k-1$ numbers $$(2\pi)^{k-1-s} L(f,s), \qquad s\in \{1,\cdots,k-1\}.$$
	 \end{example}
	 
	 \begin{remark}
	 Theorecially, it would be better to express $L$-values in terms of periods $\omega^+(f)$ and $\omega^-(f)$ of a newform. The language of $L$-values and periods are however equivalent for all cases except when weight $k=4$ and central value vanishes. For newform of level dividing $16$, which is all we will need in this article, this does not occur.
	 \end{remark}
	 
	 \begin{example}\label{k-1andk-2_CM_reduction}
	 If $f$ has CM by $\chi_D$, since $g(\chi_D) \in \mathbb{Q}\sqrt{D}$, Proposition above implies
	 $$\left(\frac{(2\pi i) L(f,k-2)}{L(f,k-1)\sqrt{D}} \right)^\sigma=  \frac{(2\pi i) L(f^\sigma,k-2)}{L(f^\sigma,k-1)\sqrt{D}},$$
	 that is, $$\frac{2\pi}{\sqrt{|D|}} \frac{L(f,k-2)}{L(f,k-1)} \in \mathbb{Q}(f).$$ The converse is also expected to hold: if $f$ is non-CM newform, then above number is transcendental.
	 \end{example}
	 
	 \begin{example}
	 The space $\Snew_7(4,\chi_{-4})$ has dimension $2$, and is spanned by two newforms $f_1,f_2$, where
	 $$f_1 = q + (\beta+2)q^2 - 4\beta q^3 + (4\beta-56)q^4 + 10q^5 + (-8\beta + 240)q^6 + 40\beta q^7 + O(q^{8}), \qquad \beta = 2\sqrt{-15},$$
	 and $f_2$ its Galois conjugate. We know, from (b) in Proposition \ref{prop_level16_twists}, $f_1 \otimes \chi_{-4} = f_2$, as can be seen by looking at first few Fourier coefficients: $$a_p(f_1) = \chi_{-4}(p) a_p(f_2),\quad p\neq 2.$$
	 Proposition \ref{Shimura_L_ratio} tells us
	 $$\frac{A(f_1,\chi_{-4},6)}{A(f_1,\mathbbm{1},5)} = \frac{L(f_2,6)}{(2\pi i) \times 2i \times L(f_1,5)} \in \mathbb{Q}(\sqrt{-15}),$$ and indeed, it equals $\frac{-597 \pm 35 \sqrt{-15}}{7808}$.
	 \end{example}
	 
	 \section{$L$-values from cuspidal subspace}
	 \subsection{Dimension counting}
	 Recall our notation of field $K$ as in equation (\ref{KandX}). We will establish in this section an upper bound for $K$-dimension of $$\text{Span}_K \left\{(2\pi)^{k-1-s} L(f,s) \Big| s\in\{1,\cdots,k-1\}, f\in \mathcal{S}_k(\Gamma;K)\right\}.$$
	 which is
	 \begin{equation}\label{dim_cusp_Lvalues}\begin{cases}2\lfloor \frac{k+2}{4}\rfloor -2 & k\text{ even and }\Gamma = \Gamma_1(4) \\ \lfloor \frac{k}{2}\rfloor - 1 & k\text{ odd and }\Gamma = \Gamma_1(4)\\ k-2 & \Gamma = \Gamma(4) \\  2k-4 & k\text{ odd and }\Gamma = \Gamma_1(8) \\ \lfloor \frac{3k}{2}\rfloor - 3& k\text{ odd and }\Gamma = \Gamma_1(8)  \end{cases}.\end{equation}
	 Since $\int_0^{i\infty} f(M\tau) d\tau = M^{-1} \int_0^{i \infty} f(\tau) d\tau$, it suffices to focus on forms in newspace. 
	 Let $\mathcal{N}$ be the set of weight $k$ newforms with level dividing that of $\Gamma$. Introduce the notation $$A'(f,\chi,s) := (2\pi)^{k-1-s} L(f\otimes \chi, s).$$
	 Let $V_{k,\Gamma}$ be the finite-dimensional $\overline{\mathbb{Q}}$-vector space spanned by formal symbols $[f,s]$, where $$s\in \{1,\cdots,k-1\}, \qquad  f\in \mathcal{N}.$$
	 Let $W_{k,\Gamma}$ be the subspace of $V_{k,\Gamma}$ spanned by 
	 $$ [f\otimes \chi, m_1] - \frac{A'(f,\chi,m_1)}{A'(f,1,m_2)} [f,m_2]$$
	 where 
	 \begin{align}\label{para_restriction}\begin{split}
	 &f\in \mathcal{N}, \chi \text{ primitive Dirichlet character such that }f\otimes \chi \in \mathcal{N};\\ 
	 &m_1,m_2 \in  \{1,\cdots,k-1\};\\ 
	 &\chi(-1) = (-1)^{m_1-m_2}, A'(f,1,m_2)\neq 0.  \end{split}
	 \end{align}
	Note that by Proposition \ref{Shimura_L_ratio}, the number $\frac{A'(f,\chi,m_1)}{A'(f,1,m_2)}$ is algebraic. The main process of dimension counting is following lemma.
	 
	 \begin{lemma}
	 The $\overline{\mathbb{Q}}$-dimension of $V_{k,\Gamma}/W_{k,\Gamma}$ is the same as equation (\ref{dim_cusp_Lvalues}). 
	 \end{lemma}
	 \begin{proof}
	 We count the number of distinct $[f,s]$ with $s\in \{1,\cdots,k-1\},  f\in \mathcal{N}$ that remains after modulo relations provided by $W_{k,\Gamma}$. Firstly, by induction on level, it suffices to focus on $f$ that are twist-minimal. Secondly, the relation in $V_{k,\Gamma}/W_{k,\Gamma}$
	 $$[f,m_1] - \frac{A'(f,1,m_1)}{A'(f,1,m_2)} [f,m_2] = 0 \qquad m_1\equiv m_2 \pmod{2}$$
	 implies that we only need to concentrate on $[f,k-1]$ and $[f,k-2]$ for each $f\in \mathcal{N}$ (Example \ref{k-1andk-2reduction}). Furthermore, if $f$ has CM by $\chi$ (which must be an odd character), then $[f,k-2]$ is a multiple of $[f,k-1]$, thus only one basis element $[f,k-1]$ remains (Example \ref{k-1andk-2_CM_reduction}). Now we proceed case by case, we will freely use statements from Lemmas \ref{prop_level16_twists} and \ref{CM_forms}.
	 \begin{enumerate}[leftmargin=0pt]
	 	\item $\Gamma = \Gamma_1(4)$ and $k$ even: in this case every newform is twist-minimal, and every 		
	 	newform gives two basis elements: $[f,k-1]$ and $[f,k-2]$. Thus the dimension of $V_{k,\Gamma}/W_{k,\Gamma}$ is
	 	$$2\sum_{d\mid 4} \dim \Snew_k(d,\mathbbm{1}) = 2\lfloor \frac{k+2}{4}\rfloor -2.$$
	 	\item $\Gamma = \Gamma_1(4)$ and $k$ odd: in this case every newform is twist-minimal. For a non-CM form $f$, $f\neq  f \otimes \chi_{-4}$, we can reduce the four basis element, $$[f,s], [f\otimes \chi_{-4},s]\qquad s\in \{k-2,k-1\},$$
	 	into two: $[f,k-1]$ and $[f\otimes \chi_{-4}, k-1]$. So the two (distinct) newforms $f$ and $f\otimes \chi_{-4}$ produce two basis elements. If $f$ has CM, then it only gives one: $[f,k-1]$. Hence the $\overline{\mathbb{Q}}$-dimension is $$\dim \Snew_k(4,\chi_{-4}) =  \lfloor \frac{k}{2}\rfloor - 1.$$
	 	\item $\Gamma = \Gamma(4)$ and $k$ even: by equation \ref{GammaNinflation}, we are essentially looking at $\mathcal{S}_k(16,\mathbbm{1})$. Those level 16 newforms are not twist-minimal, so can be ignored. Every other newform gives two basis elements: $[f,k-2]$ and $[f,k-1]$, so dimension is 
	 	$$2\sum_{d\mid 8} \dim \Snew_k(d,\mathbbm{1}) = k-2.$$
	 	\item $\Gamma = \Gamma(4)$ and $k$ odd: in this case one looks at $\mathcal{S}_k(16,\chi_{-4})$. Since $\Snew_k(8,\chi_{-4})$ is trivial, we only need to consider $f$ coming from $\Snew_k(4,\chi_{-4})$ and $\Snew_k(16,\chi_{-4})$. There is exactly one CM newform: when $k\equiv 1\pmod{4}$, it comes from $\Snew_k(4,\chi_{-4})$; when $k\equiv 3\pmod{4}$, it comes from $\Snew_k(16,\chi_{-4})$. The CM form gives one basis element. A pair of non-CM form gives two elements:  $[f,k-1]$ and $[f\otimes \chi_{-4},k-1]$. So the $\overline{\mathbb{Q}}$-dimension is $$\dim \Snew_k(4,\chi_{-4}) + \dim \Snew_k(16,\chi_{-4}) = k-2$$
	 	\item $\Gamma = \Gamma_1(8)$ and $k$ even. In this case, $f$ either lies in $\Snew_k(d,\mathbbm{1})$ for $d=1,2,4$ or $8$; or $f\in \Snew_k(8,\chi_8)$. For first possibility, the reasoning similar to that of (1) gives $2\sum_{d\mid 8} \dim \Snew_k(d,\mathbbm{1})$ different basis elements. For second possibility, we reduce the four numbers, $$[f,s], [f\otimes \chi_8,s]\qquad s\in \{k-2,k-1\},$$
	 	into two: $[f,k-1]$ and $[f,k-2]$. So the two (distinct) newforms $f$ and $f\otimes \chi_{8}$ produces two constants. Combining these two cases gives the dimension $$2\sum_{d\mid 8} \dim \Snew_k(d,\mathbbm{1})  + \dim \Snew_k(8,\chi_8) = 2k-4.$$
	 	\item $\Gamma = \Gamma_1(8)$ and $k$ odd. In this case, $f$ either lies in $\Snew_k(4,\chi_{-4})$ or $\Snew_k(8,\chi_{-8})$ because $\Snew_k(8,\chi_{-4})$ is trivial. They are stable under twisting by $\chi_{-4}$ and $\chi_{-8}$ respectively. Argument similar to that of (2) says dimension is $$\dim \Snew_k(4,\chi_{-4}) + \dim \Snew_k(8,\chi_{-8}) = \lfloor \frac{3k}{2}\rfloor - 3.$$
	 \end{enumerate}
	 \end{proof}
	 
	 There is an evaluation map $\textsf{ev}:V_{k,\Gamma} \to \mathbb{C}$, which is defined on basis element $[f,s]$ by $[f,s] \mapsto A'(f,1,s) = (2\pi)^{k-1-s} L(f,s)$, and then extended linearly to all of $V_{k,\Gamma}$. By construction, $W_{k,\Gamma}$ lies in the kernel of $\textsf{ev}$, hence induces
	 $\textsf{ev}:V_{k,\Gamma}/W_{k,\Gamma} \to \mathbb{C}$. 
	 \begin{corollary}
	 The space $$\text{Span}_{\overline{\mathbb{Q}}} \left\{(2\pi)^{k-1-s} L(f,s) \Big| s\in\{1,\cdots,k-1\}, f\in \mathcal{N}\right\}$$
	 has dimension less or equal to the number given in equation (\ref{dim_cusp_Lvalues}).
	 \end{corollary}
	 \begin{proof}
	 This is the image of $\textsf{ev}$ under $V_{k,\Gamma}$, which factors through $V_{k,\Gamma}/W_{k,\Gamma}$. Above lemma says it has equation (\ref{dim_cusp_Lvalues}) as its dimension, hence its image under $\textsf{ev}$ has at most this dimension.
	 \end{proof}

	This result is about $\overline{\mathbb{Q}}$-dimension, we now refine it to a result about $K$-dimension. 
	\subsection{Galois action on $V_{k,\Gamma}$}
	As a well-known fact about newforms, $\mathcal{N}$ is closed under Galois action, it induces an action of $G_K := \text{Gal}(\overline{\mathbb{Q}}/K)$ on $V_{k,\Gamma}$ defined by $[f,s]^\sigma = [f^\sigma,s]$. 
	\begin{lemma}
	With $K$ given by (\ref{KandX}), $W_{k,\Gamma}$ is invariant under $G_K$.
	\end{lemma}
	\begin{proof}
	We need to show for each $\sigma \in G_K$, $f,\chi,m_1,m_2$ as in equation (\ref{para_restriction}), the element 
	$$ [f\otimes \chi, m_1]^\sigma - \left(\frac{A'(f,\chi,m_1)}{A'(f,1,m_2)} \right)^\sigma [f,m_2]^\sigma =  [f^\sigma \otimes \chi^\sigma, m_1] - \left(\frac{A'(f,\chi,m_1)}{A'(f,1,m_2)} \right)^\sigma [f^\sigma ,m_2] $$
	is in $W_{k,\Gamma}$. Because
	$$ [f^\sigma \otimes \chi^\sigma, m_1] - \frac{A'(f^\sigma,\chi^\sigma,m_1)}{A'(f^\sigma,1,m_2)} [f^\sigma,m_2]$$ is in $W_{k,\Gamma}$ by definition, comparing with above expression, it suffices to show $$\left(\frac{A'(f,\chi,m_1)}{A'(f,1,m_2)} \right)^\sigma = \frac{A'(f^\sigma,\chi^\sigma,m_1)}{A'(f^\sigma,1,m_2)}.$$
	Proposition \ref{Shimura_L_ratio} implies 
	$$\left(\frac{A'(f,\chi_1,m_1)}{A'(f,1,m_2)}\right)^\sigma =\frac{A'(f^\sigma,\chi^\sigma,m_1)}{A'(f^\sigma,1,m_2)} \times (i^{m_2-m_1})^{1-\sigma} \frac{g(\chi)^\sigma}{g(\chi^\sigma)}, \qquad  \forall \sigma\in G_\mathbb{Q}.$$
	hence it suffices to show the factor $(i^{m_2-m_1})^{1-\sigma} \frac{g(\chi)^\sigma}{g(\chi^\sigma)}$ is $1$.\par
	Let $D$ be either a fundamental discriminant or $1$, the Gauss sum $g(\chi_D) \in \mathbb{Q} \sqrt{D}$, so $g(\chi_D)^\sigma/ g(\chi_D^\sigma) = \sqrt{D}^{1-\sigma}$. When $\chi = \chi_D$ is even (i.e. $D>0$), $m_1-m_2$ is even, thus $(i^{m_2-m_1})^{1-\sigma} = 1$. When $\chi = \chi_D$ is odd (i.e. $D<0$), $m_1-m_2$ is odd, $(i^{m_2-m_1})^{1-\sigma} = i^{1-\sigma}$. Thus in both cases, we have
	$$(i^{m_2-m_1})^{1-\sigma} \frac{g(\chi)^\sigma}{g(\chi^\sigma)} = \sqrt{|D|}^{1-\sigma}, \quad \forall \sigma\in G_\mathbb{Q}.$$
For $\Gamma = \Gamma_1(4)$ or $\Gamma(4)$, the only possible $D$ are $1$ and $-4$, $\sqrt{|D|} \in \mathbb{Q}$, so $\sqrt{|D|}^{1-\sigma}$ is $1$ for any $\sigma\in G_\mathbb{Q}$.
For $\Gamma = \Gamma_1(8)$, $D$ can be $1,-4,8$ or $-8$, in last two cases, $\sqrt{|D|}^{1-\sigma}$ is $1$ iff $\sigma$ fixes $\sqrt{8}$, i.e. $\sigma \in G_{\mathbb{Q}(\sqrt{2})}$.
	\end{proof}
	
	 \begin{corollary}
	 $$\text{Span}_K \left\{(2\pi)^{k-1-s} L(f,s) \Big| s\in\{1,\cdots,k-1\}, f\in \mathcal{S}_k(\Gamma;K)\right\}$$ has $K$-dimension less or equal to the number given in equation (\ref{dim_cusp_Lvalues}).
	 \end{corollary}
	 \begin{proof}
	 By above lemma, the action of $G_K$ on $V_{k,\Gamma}$ descends to the quotient $V_{k,\Gamma}/W_{k,\Gamma}$. Hence it has a basis fixed by $G_K$ (\cite[p.~40]{silverman2009arithmetic}). Let $N$ be the $\overline{\mathbb{Q}}$-span of elements in $\mathcal{N}$. The fixed space is\footnote{here we distribute the first component of $[f,s]$ linearly.} $$(V_{k,\Gamma})^{G_K} = \left\{[g,s] | s\in \{1,\cdots,k-1\}, g\in  N^{G_K}\right\},$$ $N^{G_K}$ consists of forms in $N$ with Fourier coefficient in $K$. As oldspace and newspace are invariant under Galois action, we have the claim.
	 \end{proof}
	 
	 Finally we prove our main theorems.
	 \begin{proof}[Proof of Theorem \ref{theorem_dim}]
	 From Theorem \ref{int_to_L_values}, we know
	 $$\mathcal{I}_k(\Gamma) = \text{Span}_K \left\{(2\pi)^{k-1-s} L(f,s) \Big| s\in\{1,\cdots,k-1\}, f\in \mathcal{M}_k^{(0,\infty)}(\Gamma;K)\right\}.$$
	 As Eisenstein and cusp space are both defined over $\mathbb{Q}$, this is the sum of
	 $$\text{Span}_K \left\{(2\pi)^{k-1-s} L(f,s) \Big| s\in\{1,\cdots,k-1\}, f\in \mathcal{E}_k^{(0,\infty)}(\Gamma;K) \right\}$$
	 and 
	 $$\text{Span}_K \left\{(2\pi)^{k-1-s} L(f,s) \Big| s\in\{1,\cdots,k-1\}, f\in S_k(\Gamma;K) \right\}.$$
	 The dimension coming from former space is (at most) the numbers of constants given in equation (\ref{eis_consts}),
	 and that from latter space is at most the number given by equation (\ref{dim_cusp_Lvalues}), adding up these two proves the theorem.
	 \end{proof}
	 \begin{proof}[Proof of Theorem \ref{theorem_const}]
	 The $L$-values of Dirichlet character is accounted for by equation (\ref{eis_consts}), the numbers involving $\Omega_{-4}, \Omega_{-8}$ comes from CM cusp forms (Proposition \ref{CM_L_value}).
	 \end{proof}
	 Theorem \ref{theorem_const} is much easier to prove than that of Theorem \ref{theorem_dim}, in fact a proof of former can already be obtained by examining ideas in \cite{wan2016integrals}.

\section*{Appendix: $L$-values of $\eta(2\tau)^{12}$ and a hypergeometric series}
In this short appendix, we provide evidence that the weight~6 cusp form $\eta(2\tau)^{12}$ arises from a hypergeometric motive (i.e., sharing the same Galois representation). This is noteworthy, as only very few cusp forms of weight $\geq 4$ are conjectured to be hypergeometric. For weights $\geq 7$, whether example exists remains to be seen. \par

Firstly, we establish the equality \begin{equation}\label{hyp_identity}L(\eta(2\tau)^{12},5) =\frac{\pi^6}{1024} \sum_{n\geq 0} \frac{(\frac12)_n^8}{(1)_n^8} (1+4n) = \frac{\pi^6}{1024} \: \pFq{9}{8}{\frac{1}{2},\frac{1}{2},\frac{1}{2},\frac{1}{2},\frac{1}{2},\frac{1}{2},\frac{1}{2},\frac{1}{2},\frac{5}{4}}{1,1,1,1,1,1,1,\frac{1}{4}}{1},\end{equation}
where $_9F_8$ is generalized hypergeometric series. 
\begin{proof}[Proof of (\ref{hyp_identity})]
	From Example \ref{motiv_example} above, $$L(f,5) = \frac{1}{16} \int_0^1 K(x)^2 K(1-x)^2 \:\mathrm{d}x.$$ We shall convert it into the hypergeometric series via Fourier-Legendre expansion (\cite{cantarini2019interplay}). Recall the Legendre polynomial $P_n(x) = \frac{1}{2^n n!} \frac{d^n}{dx^n} (x^2-1)^n$, it satisfies the orthogonality relation
	$$\int_0^1 P_n(2x-1) P_m(2x-1) \:\mathrm{d}x = \begin{cases}0\quad & m\neq n \\ 1/(2n+1)\quad & m=n \end{cases},$$
	with $\{P_n(2x-1)\}_{n\geq 0}$ dense in $L^2([0,1])$. For any $f\in L^2([0,1])$, we have $$f(x) = \sum_{n\geq 0} c_n P_n(2x-1) \text{ in } L^2([0,1]),\qquad c_n = (2n+1) \int_0^1 f(x) P_n(2x-1) \:\mathrm{d}x.$$
	Parseval's theorem applied to identity\footnote{which can be verified easily by creative telescoping techniques.} $$\int_0^1 K(x) K(1-x) P_n(2x-1) \:\mathrm{d}x = \begin{cases}0 \quad &n \text{ odd} \\ \frac{\pi}{8} \frac{\Gamma(1/2+n/2)^4}{\Gamma(1+n/2)^4} \quad &n\text{ even} \end{cases}$$
	gives our desired conclusion
	$$\int_0^1 K(x)^2 K(1-x)^2 \:\mathrm{d}x = \sum_{n\geq 0} (4n+1) \left( \frac{\pi}{8} \frac{\Gamma(\frac{1}{2}+n)^4}{\Gamma(1+n)^4}\right)^2.$$
\end{proof}

This proof relating $L$-value to hypergeometric series differs from the mainstream approach (as found in \cite{bonisch2023modularity, bonisch2024d}), which relies Euler integral representation and modularity of lower degree $_pF_q$. \par
Such a $L$-value identity often implies existence of corresponding super-congruence (\cite{allen2025explicit, zudilin2018hypergeometric}). Proposition above inspires us to conjecture the following numerically-supported super-congruence: let $a_n(f)$ be the $n$-th Fourier coefficient of $f = \eta(2\tau)^{12}$, then for primes $p\geq 7$, $$\sum_{n=0}^{p-1} \frac{(\frac12)_n^8}{(1)_n^8} (1+4n) \stackrel{?}{\equiv} p a_p(f) \pmod{p^6}.$$

The only other weight 6 modular form (\cite[Equation 17]{zudilin2018hypergeometric}) reported in literature that is believed to come from hypergeometric motive is $\eta(2\tau)^{12} + 32 \eta(2\tau)^4 \eta(8\tau)^8 \in \Snew_6(\Gamma_0(8)).$ Such congruences for weight $4$ cusp forms have been established in \cite{long2021supercongruences}.


\bibliographystyle{plain} 
\bibliography{../ref} 

\end{document}